\documentclass[12pt]{article}
\usepackage{amsmath, amsthm, amssymb}
\begin{document}

\begin{center}
\large{\textbf{ON SPLITTING OF EXACT DIFFERENTIAL FORMS}}

\bigskip

V.N.Dumachev \\
Voronezh institute of the MVD of Russia\\
e-mail: dumv@comch.ru
\end{center}

{\bf Abstract.} In work the internal structure of de Rham cohomology
is considered. As examples the phase flows in $\mathbb {R}^3$
admitting the Nambu poisson structure are studied.



\section{Introduction}

Let $\Lambda ^{k}(\mathbb{M})$ - be the exterior graded algebra of
differential forms with de Rham complex
\[
0\rightarrow \Lambda ^{0}(\mathbb{M})\rightarrow \Lambda ^{1}(\mathbb{M}%
)\rightarrow ...\rightarrow \Lambda ^{n-1}(\mathbb{M})\rightarrow \Lambda
^{n}(\mathbb{M})\rightarrow 0.
\]
Remember \cite{bott}, that differentual form $\omega \in \Lambda ^k(\mathbb{M})$ is called closed if $
d\omega =0$, and exact if $\omega=d\nu$ for some $\nu \in \Lambda ^{k-1}(\mathbb{M})$. The quotient of closed $k$-forms on manifolds $\mathbb{M}$ by exact $k$-forms is the $k'th$ de Rham cohomology group
\[
H^k(\mathbb{M})=\frac{\ker \left( d:\Lambda ^k(\mathbb{M})\rightarrow \Lambda ^{k+1}(\mathbb{M})\right)}
{\text{im}\left(
d:\Lambda ^{k-1}(\mathbb{M})\rightarrow \Lambda ^k(\mathbb{M})\right) }.
\]
It is known, that cohomology $\mathbb{M=R}^{n}$ is vanish. It means, that for anyone $\omega \in \Lambda ^{k}(\mathbb{R}^{n})$ such that $d\omega =0$ there exists a $\nu \in \Lambda ^{k-1}(\mathbb{R}^{n})$ such that $\omega =d\nu$. In any case the dimension of cohomology group is determined by Betti numbers:
\[
b^k=\frac{\dim \left( \omega :d\omega =0\right) }{\dim
\left( d\nu\right) }.
\]

\section{Cogomology in $\Lambda ^{2}(\mathbb{R}^{n})$}

Let's connect with de Ram complex the differential module $\{C,d\}$, then
\[
Z(C)=\ker d= \{x \in C|\; dx=0 \}\]
are called cocycles of module $\{C,d\}$ (space of closed forms),
\[B(C)=\text{im} d= dC=\{x=dy|\; y \in C\}\]
are called coboundary of module $\{C,d\}$ (space of exact forms). In given designations the group $i$-cohomology $H^i $ be the quotient $i$-cocycles by $i$-coboundary
\[H^i=Z^i/B^i.\]

Assume that $\omega\in B^2 \subset \Lambda^2(\mathbb{R}^n)$. This means that
\[
d\omega =0,\quad\omega =d\nu, \quad \text{where}\quad \nu \in \Lambda ^{1}(\mathbb{R}^{n}).
\]
The received quotient we can write as
\[
H_2^2(\mathbb{R}^n)=\frac{\ker \left( d:\Lambda ^2(\mathbb{R}
^n)\rightarrow \Lambda ^3(\mathbb{R}^n)\right) }{\text{im}\left(
d:\Lambda ^1(\mathbb{R}^n)\rightarrow \Lambda ^2(\mathbb{R}
^n)\right) }=Z^2/B^2.
\]
Note however that in $\Lambda ^2(\mathbb{R}^n)$ there exists a form
\[
\omega =\lambda _{1}\wedge \lambda _{2}, \quad \text{where}\quad \lambda _{1},\lambda _{2}\in \Lambda ^{1}(\mathbb{R}^{n})  \quad \text{such that}\quad d\lambda _i=0.
\]
Let $\lambda _i\notin H^1(\mathbb{R}^n)$ (i.e. $\lambda _i=d\mu _i)$, then
\[
d\omega =0, \qquad \omega =d\mu _1\wedge d\mu _2.
\]
In other words, exact $\omega \in\Lambda ^2(\mathbb{R}^n)$ be wedge product of the other exact forms. We can write this space as
\[
B^{1,1}=\{x=dy_1 \wedge dy_2|y_i \in C\}
\]
then
\[
H_{1,1}^{2}(\mathbb{R}^{n})=\frac{\ker \left( d:\Lambda ^{2}(\mathbb{R}%
^{n})\rightarrow \Lambda ^{3}(\mathbb{R}^{n})\right) }{\text{im}\left(
d:\Lambda ^{0}(\mathbb{R}^{n})\rightarrow \Lambda ^{1}(\mathbb{R}%
^{n})\right) \oplus \text{im}\left( d:\Lambda ^{0}(\mathbb{R}%
^{n})\rightarrow \Lambda ^{1}(\mathbb{R}^{n})\right) }=Z^2/B^{1,1},
\]
or
\[
b_{1,1}^{2}(\mathbb{R}^{n})=\frac{\dim \left( \omega :d\omega =0\right) }{%
\dim \left( d\mu \wedge d\mu \right) }.
\]
It is obvious that from $B^2=\{ d\mu_1 \wedge d\mu_2,d\nu \}$, and $B^{1,1}=\{ d\mu_1 \wedge d\mu_2 \}$ it follows that $B^{1,1} \subset B^2$. This means that quotient
\[
B^2/B^{1,1} \simeq H_{1,1}^2 / H_2^2
\]
should characterize presence of obstacles (topological defects) for existence of the exact forms in $\Lambda ^{2}(\mathbb{R}%
^{n})$,  which are wedge product of exact form from $\Lambda ^{1}(\mathbb{R}^{n})$.

\bigskip
\textbf{Example 1.} Consider the dynamical systems in $\mathbb{R}^{3}$
\[
\overset{\cdot }{x}=-xz; \qquad \overset{\cdot }{y}=yz; \qquad \overset{\cdot }{z}=x^{2}-y^{2}.
\]
According \cite{dum}, this phase flow has one vectorial
\[
\textbf{h}=\frac{1}{4}\left(
(-x^{2}y+y^{3}+yz^{2})dx+(x^{3}-y^{2}x+xz^{2})dy-2xyzdz\right)
\]
and two scalar Hamiltonians
\[
H=\frac{1}{2}(x^{2}+y^{2}+z^{2}),\qquad F=xy.
\]
These Hamiltonians are connected by expressions
\[
d\textbf{h}=dH\wedge dF.
\]
This means that our system admit Poisson structure with vectorial Hamiltonian
\[
\overset{\cdot }{x}_{i}=\{\textbf{h},x_{i}\}=X_{h}\rfloor dx_{i},
\]
where
\[
X_{h}=-xz\frac{\partial }{\partial x}+yz\frac{\partial }{\partial y}%
+(x^{2}-y^{2})\frac{\partial }{\partial z},
\]
and Poisson structure (Nambu \cite{nambu})
\[
\overset{\cdot }{x}_{i}=\{H,F,x_{i}\}=X_{H}\rfloor dF\wedge
dx_{i}=-X_{F}\rfloor dH\wedge dx_{i},
\]
where
\begin{eqnarray*}
X_{H} &=&z\frac{\partial }{\partial x}\wedge \frac{\partial }{\partial y}+x%
\frac{\partial }{\partial y}\wedge \frac{\partial }{\partial z}+y\frac{%
\partial }{\partial z}\wedge \frac{\partial }{\partial x}, \\
X_{F} &=&y\frac{\partial }{\partial y}\wedge \frac{\partial }{\partial z}+x%
\frac{\partial }{\partial z}\wedge \frac{\partial }{\partial x}.
\end{eqnarray*}

\bigskip
\textbf{Example 2.} Divergence-free Lorenz set
\[
\overset{\cdot }{x}=y-z; \qquad \overset{\cdot }{y}=-x+xz; \qquad \overset{\cdot }{z}=x-xy
\]
has one vectorial
\begin{eqnarray*}
\textbf{h}&=&\left( \frac{x}{4}(z^{2}+y^{2})-\frac{x}{3}(z+y)\right) dx\\
&+&\left( \frac{1}{
3}(x^{2}-yz+z^{2})-\frac{1}{4}x^{2}y\right) dy\\
&+&\left( \frac{1}{3}%
(y^{2}-zy+x^{2})-\frac{1}{4}x^{2}z\right) dz
\end{eqnarray*}
and two scalar Hamiltonians
\[
H=\frac{1}{2}(x^{2}+y^{2}+z^{2}),\quad F=\left(y-\frac{y^2}{2}\right)+\left(z-\frac{z^2}{2}\right)
\]
connected by expressions
\[
d\textbf{h}=dH\wedge dF.
\]
This means that our system admit Poisson structure with vectorial Hamiltonian
\[
\overset{\cdot }{x}_{i}=\{\textbf{h},x_{i}\}=X_{h}\rfloor dx_{i}
\]
where
\[
X_{h}=(y-z)\frac{\partial }{\partial x}+(-x+xz)\frac{\partial }{\partial y}%
+(x-xy)\frac{\partial }{\partial z},
\]
and Poisson structure in two forms
\[
\overset{\cdot }{x}_{i}=\{H,F,x_{i}\}=X_{H}\rfloor dF\wedge
dx_{i}=-X_{F}\rfloor dH\wedge dx_{i},
\]
where
\begin{eqnarray*}
X_{H} &=&z\frac{\partial }{\partial x}\wedge \frac{\partial }{\partial y}+x%
\frac{\partial }{\partial y}\wedge \frac{\partial }{\partial z}+y\frac{%
\partial }{\partial z}\wedge \frac{\partial }{\partial x}, \\
X_{F} &=&(1-z)\frac{\partial }{\partial x}\wedge \frac{\partial }{\partial y}+(1-y)
\frac{\partial }{\partial z}\wedge \frac{\partial }{\partial x}.
\end{eqnarray*}

\bigskip
\textbf{Example 3.} Phase flow
\[
\overset{\cdot }{x}=xy; \qquad \overset{\cdot }{y}=x-z; \qquad \overset{\cdot }{z}=-zy
\]
has one vectorial
\[
\textbf{h}=\frac{1}{12}\left(
\begin{array}{c}
z(3y^2+4x-4z \\
-6xyz \\
x(3y²+4z-4x)
\end{array}
\right),
\]
one scalar Hamiltonian
\[
H=x-\frac{y^2}{2}+z
\]
and prehamiltonian form
\[
\Theta=-zdx+xdz
\]
connected by expressions
\[
d\textbf{h}=dH\wedge \Theta.
\]
This means that our system admit Poisson structure with vectorial Hamiltonian
\[
\overset{\cdot }{x}_{i}=\{\textbf{h},x_{i}\}=X_{h}\rfloor dx_{i}
\]
where
\[
X_{h}=xy\frac{\partial }{\partial x}+(x-z)\frac{\partial }{\partial y}-yz\frac{\partial }{\partial z},
\]
and scalar Poisson structure in the form
\[
\overset{\cdot }{x}_{i}=X_{H}\rfloor \Theta\wedge
dx_{i},
\]
where
\begin{eqnarray*}
X_{H} &=&\frac{\partial }{\partial x}\wedge \frac{\partial }{\partial y}+
\frac{\partial }{\partial y}\wedge \frac{\partial }{\partial z}-y\frac{
\partial }{\partial z}\wedge \frac{\partial }{\partial x}.
\end{eqnarray*}
Therefore,
\begin{equation*}
dh\in B^2 \setminus B^{1,1}.
\end{equation*}
For completeness of a statement we shall notice that as $d\Theta\neq 0$, but $d\Theta \wedge \Theta =0$, then Pfaff equation on prehamiltonians form $\Theta$ has solved with integrating factor
\[
dF=\frac{\Theta}{x^2 +z^2}, \qquad \Rightarrow \qquad F=\arctan \frac{z}{x}.
\]
The caused of global non-integrability of the given system is the presence holes $(x=0,\;z=0)$ in $x0z$ planes. It is obvious that
vanishing of second Hamiltonians has not admitted to enter of Nambu structure with a bracket $\{H,F,G\}$.

\section{Cogomology in $\Lambda ^{3}(\mathbb{R}^{n})$}

Further, the top index of the any form will denote its degree, i.e.
 $\omega ^{k}\in \Lambda ^{k}(\mathbb{R}^{n})$. For standard de Rham complexes in $\mathbb{R%
}^{n}$\ we get $\omega ^{3}\in B^3 \subset
\Lambda ^{3}(\mathbb{R}^{n})$. This means that
\[
d\omega ^{3}=0,\quad\omega ^{3}=d\nu ^{2},\quad \text{where}\quad \nu ^{2}\in \Lambda ^{2}(\mathbb{R}^{n}).
\]
Thus
\begin{eqnarray*}
H_{3}^{3}(\mathbb{R}^{n})&=&\frac{\ker \left( d:\Lambda ^{3}(\mathbb{R}%
^{n})\rightarrow \Lambda ^{4}(\mathbb{R}^{n})\right) }{\text{im}\left(
d:\Lambda ^{2}(\mathbb{R}^{n})\rightarrow \Lambda ^{3}(\mathbb{R}%
^{n})\right) }=Z^3/B^3,\\
\\
b_3^3&=&\frac{\dim \left( \omega ^{3}:d\omega ^{3}=0\right) }{\dim
\left( d\nu ^{2}\right) }.
\end{eqnarray*}
But in $\Lambda ^{3}(\mathbb{R}^{n})$ there exists are forms
\[
\omega _{1}^{3}=\lambda _{1}^{1}\wedge \lambda _{2}^{1}\wedge \lambda
_{3}^{1}, \quad \text{where} \quad \lambda _{i}^{1}\in \Lambda
^{1}(\mathbb{R}^{n}), \quad \text{such that} \quad d\lambda
_{i}^{1}=0.
\]
Let $\lambda _{i}^{1}\notin H^{1}(\mathbb{R}^{n})$ (i.e. $\lambda _{i}^{1}=d\mu _{i}^{0})$, then
\[
d\omega _{1}^{3} =0, \qquad \omega _{1}^{3} =d\mu _{1}^{0}\wedge d\mu _{2}^{0}\wedge d\mu _{3}^{0}.
\]
In other words, exact $\omega _{1}^{3}\in \Lambda ^{3}(\mathbb{R}^{n})$  be wedge product of the other exact forms. We can write this quotient as
\[
H_{1,1,1}^{3}(\mathbb{R}^{n})=\frac{\ker \left( d:\Lambda ^{3}(\mathbb{R}%
^{n})\rightarrow \Lambda ^{4}(\mathbb{R}^{n})\right) }{\bigoplus%
\limits_{k=1}^{3}\text{im}_{k}\left( d:\Lambda ^{0}(\mathbb{R}%
^{n})\rightarrow \Lambda ^{1}(\mathbb{R}^{n})\right) }=Z^3/B^{1,1,1},
\]
\[
b_{1,1,1}^{3}=\frac{\dim \left( \omega ^{3}:d\omega
^{3}=0\right) }{\dim \left( d\mu ^{0}\wedge d\mu ^{0}\wedge d\mu ^{0}\right)
}.
\]
At the same time in $\Lambda ^{3}(\mathbb{R}
^{n})$ it is possible also to construct the forms
\[
\omega _{2}^{3}=\lambda ^{1}\wedge \lambda ^{2},\quad \text{where}\quad \lambda ^{k}\in \Lambda ^{k}(\mathbb{R}^{n})
\]
such that $d\lambda ^{k}=0$. Let $\lambda ^{k}\notin H^{k}(\mathbb{R}^{n})$ (i.e. $\lambda _{i}^{k}=d\mu _{i}^{k})$, then
\[
d\omega _{2}^{3} =0, \qquad \omega _{2}^{3} =d\mu ^{0}\wedge d\mu ^{1}.
\]
In other words, exact $\omega _{2}^{3}$  be wedge product of the other exact forms. We can write this quotient as
\[
H_{2,1}^{3}(\mathbb{R}^{n})=\frac{\ker \left( d:\Lambda ^{3}(\mathbb{R}%
^{n})\rightarrow \Lambda ^{4}(\mathbb{R}^{n})\right) }{\text{im}\left(
d:\Lambda ^{0}(\mathbb{R}^{n})\rightarrow \Lambda ^{1}(\mathbb{R}%
^{n})\right) \oplus \text{im}\left( d:\Lambda ^{1}(\mathbb{R}%
^{n})\rightarrow \Lambda ^{2}(\mathbb{R}^{n})\right) }=Z^3/B^{2,1},
\]
\[
b_{1,2}^{3}=\frac{\dim \left(
\omega ^{3}:d\omega ^{3}=0\right) }{\dim \left( d\mu ^{0}\wedge d\mu
^{1}\right) }.
\]
Evidently that
\[
H_{1,1,1}^{3}(\mathbb{R}^{n})\supset H_{1,2}^{3}(\mathbb{R}^{n})\supset
H_{3}^{3}(\mathbb{R}^{n}),
\]
such that quotients $H_{1,1,1}^{3}(\mathbb{R}^{n})/H_{3}^{3}(\mathbb{R}^{n})$ and
$H_{1,1,1}^{3}(\mathbb{R}^{n})/H_{2,1}^{3}(\mathbb{R}^{n})$
should characterize presence of obstacles (topological defects) for existence of the exact forms in
 $\Lambda ^{3}(\mathbb{R}^{n})$, which are wedge product of exact form from $
\Lambda ^{1}(\mathbb{R}^{n})$ or from $\Lambda ^{2}(\mathbb{R}^{n})$.

\section{Cogomology in $\Lambda ^{k}(\mathbb{R}^{n})$}

Generalizing the previous calculations we consider de Rham complex in $\mathbb{R}^{n}$ and get $\omega ^{k}\in \Lambda ^{k}(\mathbb{R}^{n})\notin H^{k}(\mathbb{R}^{n})$.
Then
\[
d\omega ^{k}=0,\quad\omega ^{k}=d\nu ^{k-1},\quad \text{where}\quad \nu ^{k-1}\in \Lambda ^{k-1}(\mathbb{R}^{n}),
\]
and
\[
H_{k}^{k}(\mathbb{R}^{n})=\frac{\ker \left( d:\Lambda ^{k}(\mathbb{R}%
^{n})\rightarrow \Lambda ^{k+1}(\mathbb{R}^{n})\right) }{\text{im}\left(
d:\Lambda ^{k-1}(\mathbb{R}^{n})\rightarrow \Lambda ^{k}(\mathbb{R}%
^{n})\right) }=Z^k/B^k.
\]
However, into $\Lambda ^{k}(\mathbb{R}^{n})$ there exists are forms
\[
\omega ^{k}=\bigwedge_{i=1}^{k}\lambda _{i}^{1},\quad \text{where}\quad \lambda ^{1}\in \Lambda ^{1}(\mathbb{R}^{n})
\]
such that $d\lambda ^{1}=0$. Suppose that $\lambda ^{1}\notin H^{1}(\mathbb{R}^{n})$ (i.e. $\lambda ^{1}=d\mu ^{0})$, then
\[
d\omega ^{k} =0, \qquad \omega ^{k} =\bigwedge_{i=1}^{k}d\mu^0_i.
\]
In other words, exact $\omega ^{k}\in \Lambda ^{k}(\mathbb{R}^{n})$ be wedge product of the other exact forms. We can write this quotient as
\[
H_{1,1,...,1}^{k}(\mathbb{R}^{n})=\frac{\ker \left( d:\Lambda ^{k}(\mathbb{R}%
^{n})\rightarrow \Lambda ^{n+1}(\mathbb{R}^{n})\right) }{\bigoplus%
\limits_{i=1}^{k}\text{im}_{i}\left( d:\Lambda ^{0}(\mathbb{R}%
^{n})\rightarrow \Lambda ^{1}(\mathbb{R}^{n})\right) }=Z^k/B^{1,1,...,1},
\]
\[
b_{1,1,...,1}^{k}=\frac{\dim \left(
\omega ^{k}:d\omega ^{k}=0\right) }{\dim \left( \bigwedge_{i=1}^{k}d\mu
_{i}\right) }.
\]

Continuing it is similarly we receive:

\noindent for $\omega ^{k}=d\mu ^1\wedge \bigwedge_{i=1}^{k-2}d\mu
_{i}^{0}$%
\[
H_{2,1,1,...,1}^{k}(\mathbb{R}^{n})=\frac{\ker \left( d:\Lambda ^{k}(\mathbb{%
R}^{n})\rightarrow \Lambda ^{k+1}(\mathbb{R}^{n})\right) }{\text{im}\left(
d:\Lambda ^{1}(\mathbb{R}^{n})\rightarrow \Lambda ^{2}(\mathbb{R}%
^{n})\right) \bigoplus\limits_{i=1}^{k-2}\text{im}_{i}\left( d:\Lambda ^{0}(%
\mathbb{R}^{n})\rightarrow \Lambda ^{1}(\mathbb{R}^{n})\right) },
\]
\[
b_{2,1,1,...,1}^{k}=\frac{\dim
\left( \omega ^{k}:d\omega ^{k}=0\right) }{\dim \left( d\mu ^{1}\wedge
\bigwedge_{i=1}^{k-2}d\mu _{i}^{0}\right) },
\]
for $\omega ^{k}=d\mu ^{2}\wedge \bigwedge_{i=1}^{k-3}d\mu
_{i}^{0}$%
\[
H_{3,1,1,...,1}^{k}(\mathbb{R}^{n})=\frac{\ker \left( d:\Lambda ^{k}(\mathbb{%
R}^{n})\rightarrow \Lambda ^{k+1}(\mathbb{R}^{n})\right) }{\text{im}\left(
d:\Lambda ^{2}(\mathbb{R}^{n})\rightarrow \Lambda ^{3}(\mathbb{R}%
^{n})\right) \bigoplus\limits_{i=1}^{k-3}\text{im}_{i}\left( d:\Lambda ^{0}(%
\mathbb{R}^{n})\rightarrow \Lambda ^{1}(\mathbb{R}^{n})\right) },
\]
where
\[
b_{3,1,1,...,1}^{k}=\frac{\dim
\left( \omega ^{k}:d\omega ^{k}=0\right) }{\dim \left( d\mu ^{2}\wedge
\bigwedge_{i=1}^{k-3}d\mu _{i}^{0}\right) };
\]
for $\omega ^{k}=d\mu ^{3}\wedge \bigwedge_{i=1}^{k-4}d\mu
_{i}^{0}$%
\[
H_{4,1,1,...,1}^{k}(\mathbb{R}^{n})=\frac{\ker \left( d:\Lambda ^{k}(\mathbb{%
R}^{n})\rightarrow \Lambda ^{k+1}(\mathbb{R}^{n})\right) }{\text{im}\left(
d:\Lambda ^{3}(\mathbb{R}^{n})\rightarrow \Lambda ^{4}(\mathbb{R}%
^{n})\right) \bigoplus\limits_{i=1}^{k-4}\text{im}_{i}\left( d:\Lambda ^{0}(%
\mathbb{R}^{n})\rightarrow \Lambda ^{1}(\mathbb{R}^{n})\right) },
\]
where
\[
b_{4,1,1,...,1}^{k}=\frac{\dim
\left( \omega ^{k}:d\omega ^{k}=0\right) }{\dim \left( d\mu ^{3}\wedge
\bigwedge_{i=1}^{k-4}d\mu _{i}^{0}\right) };
\]

\begin{center}
...   ...   ...
\end{center}

for $\omega ^{k}=d\mu ^{1}\wedge d\mu ^{1}\wedge
\bigwedge_{i=1}^{k-4}d\mu _{i}^{0}$%
\[
H_{2,2,1,...,1}^{k}(\mathbb{R}^{n})=\frac{\ker \left( d:\Lambda ^{k}(\mathbb{%
R}^{n})\rightarrow \Lambda ^{k+1}(\mathbb{R}^{n})\right) }{%
\bigoplus\limits_{i=1}^{2}\text{im}_{i}\left( d:\Lambda ^{1}(\mathbb{R}%
^{n})\rightarrow \Lambda ^{2}(\mathbb{R}^{n})\right)
\bigoplus\limits_{i=1}^{k-4}\text{im}_{i}\left( d:\Lambda ^{0}(\mathbb{R}%
^{n})\rightarrow \Lambda ^{1}(\mathbb{R}^{n})\right) },
\]
where
\[
b_{2,2,1,1,...,1}^{k}=\frac{\dim \left(
\omega ^{k}:d\omega ^{k}=0\right) }{\dim \left( d\mu ^{1}\wedge d\mu
^{1}\wedge \bigwedge_{i=1}^{k-4}d\mu _{i}^{0}\right) }.
\]
\begin{center}
...   ...   ...
\end{center}
etc.

For simplification of record we shall enter a multiindex
\[
\#m=\{m_{1},m_{2},...,m_{i}\},\quad m_{1}+m_{2}+...+m_{i}=k,\quad m_{1}\geq m_{2}\geq ...\geq m_{i}\geq 0,
\]
which is formed by a rule of construction of the Young diagrams. Then
\[
H_{\#m}^k(\mathbb{R}^{n})=\frac{\ker \left( d:\Lambda ^{k}(\mathbb{R}%
^{n})\rightarrow \Lambda ^{k+1}(\mathbb{R}^{n})\right) }{\bigoplus\limits_{%
\#m}\text{im}_{\#m}\left( d:\Lambda ^{\#m}(\mathbb{R}^{n})\rightarrow
\Lambda ^{\#m+1}(\mathbb{R}^{n})\right) }=Z^k/B^{\#m},
\]
\[
b_{\#m}^k=\frac{\dim \left(
\omega ^{k}:d\omega ^{k}=0\right) }{\dim \left( \bigwedge^{\#m}d\mu _{i}^{0}\right) }.
\]

So, for $k=3 $ we shall receive

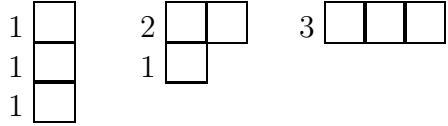
\begin{picture}(40,60)
\put(100,32){$1$}
\put(100,17){$1$}
\put(100,2){$1$}
\put(110,30){\framebox(15,15)}
\put(110,15){\framebox(15,15)}
\put(110,0 ){\framebox(15,15)}

\put(150,32){$2$}
\put(150,17){$1$}
\put(160,30){\framebox(15,15)}
\put(160,15){\framebox(15,15)}
\put(175,30 ){\framebox(15,15)}

\put(210,32){$3$}
\put(220,30){\framebox(15,15)}
\put(235,30){\framebox(15,15)}
\put(250,30 ){\framebox(15,15)}
\end{picture}

i.e.
\begin{eqnarray*}
\#m&=&\{1,1,1\} \quad \text{or} \quad H_{\#m}^k(\mathbb{R}^{n})=H_{1,1,1}^{3}(\mathbb{R}^{n}),\\
\#m&=&\{2,1\} \quad \quad \text{or} \quad H_{\#m}^k(\mathbb{R}^{n})=H_{2,1}^{3}(\mathbb{R}^{n}),\\
\#m&=&\{3\} \quad\quad \quad  \text{or} \quad H_{\#m}^k(\mathbb{R}^{n})=H_3^3(\mathbb{R}^{n}).
\end{eqnarray*}
Then the filtered complex of cohomology can be represented as follows
\begin{equation*}
H_{1,1,1}^{3}\rightarrow H_{2,1}^{3}\rightarrow H_{3}^{3}.
\end{equation*}

For $k=4$ we get

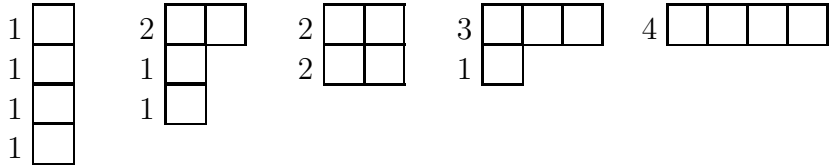
\begin{picture}(40,70)
\put(10,47){$1$}
\put(10,32){$1$}
\put(10,17){$1$}
\put(10,2){$1$}
\put(20,45){\framebox(15,15)}
\put(20,30){\framebox(15,15)}
\put(20,15){\framebox(15,15)}
\put(20,0 ){\framebox(15,15)}

\put(60,47){$2$}
\put(60,32){$1$}
\put(60,17){$1$}
\put(70,45){\framebox(15,15)}
\put(85,45 ){\framebox(15,15)}
\put(70,30){\framebox(15,15)}
\put(70,15 ){\framebox(15,15)}

\put(120,47){$2$}
\put(120,32){$2$}
\put(130,45){\framebox(15,15)}
\put(130,30){\framebox(15,15)}
\put(145,45){\framebox(15,15)}
\put(145,30){\framebox(15,15)}

\put(180,47){$3$}
\put(180,32){$1$}
\put(190,45){\framebox(15,15)}
\put(205,45){\framebox(15,15)}
\put(220,45){\framebox(15,15)}
\put(190,30){\framebox(15,15)}

\put(250,47){$4$}
\put(260,45){\framebox(15,15)}
\put(275,45){\framebox(15,15)}
\put(290,45){\framebox(15,15)}
\put(305,45){\framebox(15,15)}

\end{picture}

\noindent i.e. $\#m=\{1,1,1,1\}$, $\#m=\{2,1,1\}$, $\#m=\{2,2\}$, $\#m=\{3,1\}$ or $\#m=\{4\}$. Let's notice, that at $k\geq 4 $ the structure of cohomology $H_{\#m}^k(\mathbb{R}^{n})$ is not linear and for some small $k=4,5,6$ is shown in figures:

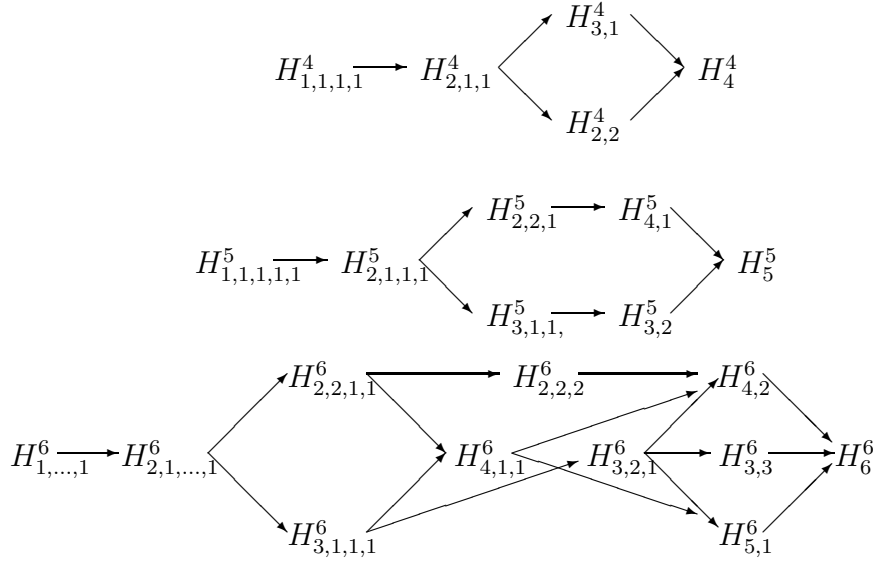
\begin{picture}(40,60)
\put(100,20){$H^4_{1,1,1,1} $}
\put(130,25){\vector(1,0){20}} \put(155,20){$H^4_{2,1,1} $}
\put(185,25){\vector(1,1){20}}\put(185,25){\vector(1,-1){20}}
\put(210,40){$H^4_{3,1} $}\put(235,45){\vector(1,-1){20}}
\put(210,0) {$H^4_{2,2} $}\put(235,5){\vector(1,1){20}}
\put(260,20){$H^4_{4} $}
\end{picture}

\bigskip
\begin{picture}(40,60)
\put(70,20){$H^5_{1,1,1,1,1} $}\put(100,25){\vector(1,0){20}} \put(125,20){$H^5_{2,1,1,1} $}
\put(155,25){\vector(1,1){20}}\put(155,25){\vector(1,-1){20}}
\put(180,40){$H^5_{2,2,1} $}\put(205,45){\vector(1,0){20}}\put(230,40){$H^5_{4,1} $}\put(250,45){\vector(1,-1){20}}

\put(180,0) {$H^5_{3,1,1,} $}\put(205,5){\vector(1,0){20}}\put(230,0) {$H^5_{3,2} $}\put(250,5){\vector(1,1){20}}

\put(275,20){$H^5_5 $}
\end{picture}

\bigskip
\begin{picture}(30,70)
\put(0,30){$H^6_{1,...,1} $}\put(18,35){\vector(1,0){23}} \put(42,30){$H^6_{2,1,...,1} $}
\put(75,35){\vector(1,1){30}}\put(75,35){\vector(1,-1){30}}

\put(105,60){$H^6_{2,2,1,1} $} \put(135,65){\vector(1,0){50}}\put(190,60){$H^6_{2,2,2} $}
\put(215,65){\vector(1,0){45}}
\put(135,65){\vector(1,-1){30}}\put(267,60) {$H^6_{4,2} $}\put(285,65){\vector(1,-1){26}}

\put(168,30) {$H^6_{4,1,1}$}
\put(190,35){\vector(3,-1){70}}\put(190,35){\vector(3,1){70}}
\put(218,30){$H^6_{3,2,1}$}
\put(240,35){\vector(1,-1){28}}\put(240,35){\vector(1,1){28}}\put(240,35){\vector(1,0){25}}
\put(268,30) {$H^6_{3,3} $}\put(287,35){\vector(1,0){25}}\put(312,30) {$H^6_{6} $}

\put(105,0) {$H^6_{3,1,1,1} $}\put(135,5){\vector(1,1){30}}
\put(135,5){\vector(3,1){80}}\put(268,0) {$H^6_{5,1} $}\put(285,5){\vector(1,1){26}}
\end{picture}

\bigskip
We can see that in the general case the cohomology sequences are not filtered. This means to define cohomology of cogomology is obviously impossible.


\begin{thebibliography}{99}

\bibitem{bott} R. Bott, L.W. Tu, {\it Differential Forms in Algebraic Topology. Graduate Texts in Mathematics 82}, Springer-Verlag, New York (1982).

\bibitem{dum} V.N. Dumachev, Phase flows and vector Lagrangians in $J^3(\pi)$, {\it Int.J.Pure and Appl.Math.}, {\bf 55} (2009),
147-152. (math.DG/1010.0081)

\bibitem{nambu} Y. Nambu, Generalized Hamiltonian dynamics, {\it Physical Review D.}, {\bf 7} (1973), 2405-2412.



\end{thebibliography}
\end{document}